# INFERENCE FOR MIXTURES OF SYMMETRIC DISTRIBUTIONS

By David R. Hunter, Shaoli Wang and Thomas P.
Hettmansperger

*Pennsylvania State University*

This article discusses the problem of estimation of parameters in finite mixtures when the mixture components are assumed to be symmetric and to come from the same location family. We refer to these mixtures as *semi-parametric* because no additional assumptions other than symmetry are made regarding the parametric form of the component distributions. Because the class of symmetric distributions is so broad, identifiability of parameters is a major issue in these mixtures. We develop a notion of identifiability of finite mixture models, which we call *k-identifiability*, where $k$ denotes the number of components in the mixture. We give sufficient conditions for $k$-identifiability of location mixtures of symmetric components when $k = 2$ or 3. We propose a novel distance-based method for estimating the (location and mixing) parameters from a $k$-identifiable model and establish the strong consistency and asymptotic normality of the estimator. In the specific case of $L_2$-distance, we show that our estimator generalizes the Hodges–Lehmann estimator. We discuss the numerical implementation of these procedures, along with an empirical estimate of the component distribution, in the two-component case. In comparisons with maximum likelihood estimation assuming normal components, our method produces somewhat higher standard error estimates in the case where the components are truly normal, but dramatically outperforms the normal method when the components are heavy-tailed.

**1. Introduction.** Given a random sample $X_1, \ldots, X_n$ from a symmetric distribution, Hodges and Lehmann [8] proposed an estimator for the center of symmetry, $\mu$, that consists of the median of all $n + \binom{n}{2}$ pairwise means $(X_i + X_j)/2$ for $i \leq j$. By the well-known property that a sample median minimizes the sum of absolute deviations from all the points in the sample, we may express this Hodges–Lehmann estimator as

$$\hat{\mu}_{\mathrm{HL}} = \arg\min_{\mu} \sum_{i \leq j} \left| \frac{X_i + X_j}{2} - \mu \right|. \tag{1}$$

This article extends the idea of Hodges and Lehmann to a more general setting in which the sample $X_1, \ldots, X_n$ comes not from a symmetric distribution, but from a finite mixture of location-shifted symmetric distributions. Yet, we do far more in this article than generalize the Hodges–Lehmann estimator. We propose a general method of estimation for location mixtures of symmetric components and discuss the central issue of identifiability, without which the very concept of estimation in these models is ill-defined.

As motivating examples, consider the two samples depicted in Figure 1. In the Old Faithful dataset, measurements give time in minutes between eruptions of the Old Faithful geyser in Yellowstone National Park, USA. These data are included as part of the R and S-PLUS statistics packages [type "help(faithful)" in R or "help(geyser)" in S-PLUS for more details]. For the Old Faithful eruption data, a two-component mixture model is clearly reasonable. However, in the case of the double exponential dataset, this fact is not so clear. These data were simulated from a 2-component mixture of location-shifted double exponential distributions. A common choice for fitting parameters in a 2-component mixture when nothing is known about the shape of the component distributions is to use maximum likelihood estimation based on normally distributed components. For the Old Faithful data, the method we propose in this article performs nearly identically to the normal method; furthermore, our method enables a validation of the normality assumption by providing a nonparametric estimate of the distribution function of the mixture components. For the simulated double exponential mixture on the right in Figure 1, estimates of $\mu_1 = -1$, $\mu_2 = 1$, and $\lambda_1 = 0.3$ are $(-1.04, 0.97, 0.33)$ for our method and $(-5.48, 0.33, 0.006)$ for the normal method. In these examples, our method complements the likelihood-based normal methods when they appear to be appropriate and it outperforms the normal methods when they are not appropriate. We explore this comparison in Sections 5 and 6.

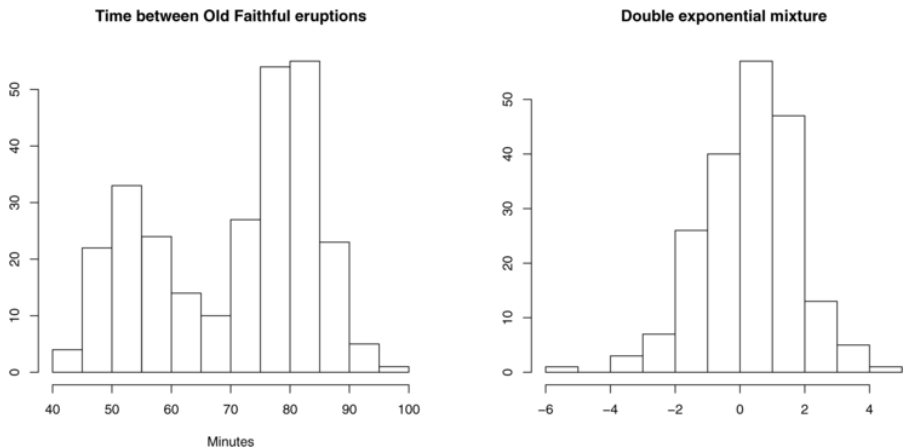

FIG. 1.   *The Old Faithful dataset is clearly suggestive of a 2-component mixture of symmetric components. The data on the right are simulated from a 2-component mixture of double exponential distributions with centers $\mu_1 = -1$ and $\mu_2 = 1$ and mixing parameters $\lambda_1 = 0.3$ and $\lambda_2 = 0.7$.*

But before we discuss the application of our method to data, there is much preparation to be done. As a first step, we formally specify the model for the data in Section 2. We also discuss the all-important topic of identifiability in Section 2. Parameter estimation is the topic of Section 3, where we propose a class of estimators and prove that they are strongly consistent. Section 4 explores the connection between a particular form of our estimator and the Hodges–Lehmann estimator, exploiting this connection to establish asymptotic normality. Section 5 examines the numerical implementation of our estimation method, giving several examples. Finally, the discussion in Section 6 compares our estimation method to the canonical estimation method for problems of this type, namely, maximum likelihood estimation assuming a mixture of normal distributions. Technical details about identifiability and proofs of strong consistency are given in Appendices A and B, respectively.

**2. The model and identifiability.** Suppose that $X_1, \ldots, X_n$ are independent and identically distributed from a $k$-component mixture distribution with distribution function

$$(2) \qquad F(x) = \sum_{j=1}^{k} \lambda_j G(x - \mu_j)$$

for some distribution function $G(x)$ that is completely unspecified, except for the assumption that $G$ is symmetric about zero, that is $G(x) = 1 - G(-x)$ for all continuity points $x$ of $G$. In this article, we denote by $\mathcal{S}$ the set of all distributions symmetric about zero and we refer to such distributions as *zero-symmetric*. The shifted distributions $\{G(x - \mu) : G \in \mathcal{S}, \mu \in R\}$ are referred to as *symmetric* distributions. We assume throughout that $k$ is fixed and known. Such an assumption is often justified on the basis of theory specific to the application at hand; see, for example, [7].

The distribution of equation (2) may be written as the convolution of $G$ with a distribution supported on the $k$ points $(\mu_1, \ldots, \mu_k)$. Because such finite distributions arise frequently in this article, we introduce a notation for them. For $\boldsymbol{\lambda} = (\lambda_1, \ldots, \lambda_k)$ and $\boldsymbol{\mu} = (\mu_1, \ldots, \mu_k)$ such that $\lambda_j \geq 0$ for all

$j$ and $\sum_j \lambda_j = 1$, let

$$\Delta_k(x; \boldsymbol{\lambda}, \boldsymbol{\mu}) = \sum_{j=1}^k \lambda_j \delta_{\mu_j}(x),$$

where $\delta_t(x) = I\{t \le x\}$ denotes the distribution function which assigns mass 1 to the point $t$. We sometimes abuse notation and refer to the distribution $\Delta_k(\boldsymbol{\lambda}, \boldsymbol{\mu})$ without its argument $x$. Thus, we may rewrite equation (2) as

$$(3) \qquad\qquad F = G \star \Delta_k(\boldsymbol{\lambda}, \boldsymbol{\mu}).$$

Throughout this article, we use the convention that a distribution function superscripted with a minus sign denotes the result of that distribution being reflected over the origin. For example, $\Delta_k^-(\boldsymbol{\lambda}, \boldsymbol{\mu})$ denotes the same distribution as $\Delta_k(\boldsymbol{\lambda}, -\boldsymbol{\mu})$.

Compared to the large statistical literature regarding mixture models in which the component distributions are assumed to come from a known parametric family, relatively little work has been done in the case where minimal assumptions are made regarding $G$. Hettmansperger and Thomas [7] report promising results when the component distributions are multivariate, by reducing the model to a mixture of multinomials using cutpoints. In their case, the data consist of vectors of repeated measures, assumed to be independent and identically distributed, conditional on the component from which they are drawn. Hall and Zhou [6] discuss a related situation in which the repeated measures are independent but not identically distributed and the mixture has two components. Identifiability issues make both of these approaches impossible in the univariate (nonrepeated measures) case.

Here, we take a qualitatively different approach. We consider univariate, rather than multivariate, data and we achieve identifiability by imposing a symmetry restriction on the individual components; see [3] for an alternative approach to the same problem. Cruz-Medina and Hettmansperger [4] apply the cutpoint approach of Hettmansperger and Thomas [7] to a similar case in which the component distributions are unimodal and continuous (conditions we do not assume here) in addition to being symmetric. Ellis [5] considers the problem of deconvolving $F = G \star Q$ into a symmetric part $G$ and a nonsymmetric part $Q$, but without the assumption that $Q$ is a $k$-point distribution. Finally, Walther [14, 15] considers the problem not of estimation, but rather detection of the presence of mixing for a univariate distribution under very minimal assumptions on the component distributions, although these assumptions are quite different from the assumptions we make here.

The issue of identifiability looms large in the study of mixture models—in order for parameter estimation to make any sense, we must be assured that the parameters are uniquely determined by the mixture. Here, "the mixture" is $F(x)$ from equation (3). Clearly, a permutation applied to the entries of $\boldsymbol{\lambda}$ and $\boldsymbol{\mu}$ does not change $F(x)$, but this particular identifiability conundrum—often called the "label-switching" problem—is easily solved in this case by insisting that $\mu_1 < \cdots < \mu_k$. Thus, we define the parameter space of interest to be $\Omega_k \times \mathcal{S}$, where

$$\Omega_k = \left\{ (\lambda_1, \ldots, \lambda_k; \mu_1, \ldots, \mu_k) : \lambda_j \geq 0, \sum_j \lambda_j = 1, \mu_1 < \cdots < \mu_k \right\}$$

and $\mathcal{S}$ is the set of all zero-symmetric probability distributions on the real numbers. Furthermore, let

$$\mathcal{M}_k = \left\{ F : F(x) = \sum_j \lambda_j G(x - \mu_j), (\boldsymbol{\lambda}, \boldsymbol{\mu}) \in \Omega_k, G \in \mathcal{S} \right\}$$

be the set of all mixture distributions defined by equation (3).

Typically, "identifiability" is a property of the whole set of distributions defined by a mixture model [10, 11, 13, 16]. Thus, to adopt the traditional view is to view identifiability of $\mathcal{M}_k$ as an "all-or-nothing" proposition: either $\mathcal{M}_k$ is identifiable, or it is not. (From this perspective, for $k > 1$, it is not.) However, we prefer to define identifiability as a property of individual distributions in $\mathcal{M}_k$. This view makes it possible to refer to subsets of identifiable mixture distributions within $\mathcal{M}_k$.

To develop this notion of identifiability, let $\varphi_k : \Omega_k \times \mathcal{S} \to \mathcal{M}_k$ denote the function that maps $(\boldsymbol{\lambda}, \boldsymbol{\mu}, G)$ onto $G \star \Delta_k(\boldsymbol{\lambda}, \boldsymbol{\mu})$. Essentially, identifiability means that $\varphi_k$ should be a one-to-one (i.e., an invertible) function. Let

$$\varphi_k^{-1}(F) = \{ (\boldsymbol{\lambda}, \boldsymbol{\mu}, G) \in \Omega_k \times \mathcal{S} : \varphi_k(\boldsymbol{\lambda}, \boldsymbol{\mu}, G) = F \}$$

denote the inverse image of $F$ under $\varphi_k(\cdot)$. Although $\varphi_k^{-1}(F)$ is not always a singleton for $F \in \mathcal{M}_k$, there are elements $F \in \mathcal{M}_k$ for which $\varphi_k^{-1}(F)$ is a singleton and these are precisely the $k$-component mixtures we consider $k$-identifiable.

DEFINITION 1. If $F \in \mathcal{M}_k$ has the property that $\varphi_k^{-1}(F)$ contains a single element of $\Omega_k \times \mathcal{S}$, then $F$ is said to be *identifiable as a $k$-component mixture of distributions from a symmetric location family*. Alternatively, we say that such an $F$ is *$k$-identifiable*.

In estimation, the primary interest is typically in the values of $\boldsymbol{\lambda}$ and $\boldsymbol{\mu}$. Thus, we turn our attention to the largest subset $\Omega_k^* \subset \Omega_k$ such that the image of $\Omega_k^* \times \mathcal{S}$ under the map $\varphi_k$ consists entirely of $k$-identifiable distributions,

(4)  $\Omega_k^* = \{ (\boldsymbol{\lambda}, \boldsymbol{\mu}) \in \Omega_k : G \star \Delta_k(\boldsymbol{\lambda}, \boldsymbol{\mu})$ is $k$-identifiable for all $G \in \mathcal{S} \}.$

Note that $\Omega_k^*$ is a proper subset of $\Omega_k$ for $k \geq 2$ because no element of $\Omega_k$ for which some $\lambda_j = 0$ can be in $\Omega_k^*$. By the same reasoning, a $k$-component model can always be made into a $(k + \ell)$-component model, $\ell > 0$, by adding $\ell$ components with zero weight. Therefore, no distribution can be $k$-identifiable for more than one value of $k$.

Even if all $\lambda_j$ are nonzero, a distribution is not necessarily $k$-identifiable. As an example, let $G_1(t) = \frac{1}{2}\delta_{-1}(t) + \frac{1}{2}\delta_1(t)$ be the zero-symmetric distribution with jumps of $\frac{1}{2}$ at $-1$ and $1$. We see that $G_1(t-1)$ assigns equal mass to the points $0$ and $2$, whereas $G_1(t-5)$ assigns equal mass to $4$ and $6$. Now, let $G_2(t) = \frac{1}{2}\delta_{-2}(t) + \frac{1}{2}\delta_2(t)$ be the zero-symmetric distribution with jumps of $\frac{1}{2}$ at $-2$ and $2$. This implies that $G_2(t-2)$ assigns equal mass to the points $0$ and $4$, whereas $G_2(t-4)$ assigns equal mass to $2$ and $6$. We conclude that the mixtures $\frac{1}{2}G_1(t-1) + \frac{1}{2}G_1(t-5)$ and $\frac{1}{2}G_2(t-2) + \frac{1}{2}G_2(t-4)$ both assigns mass $\frac{1}{4}$ to the points $0$, $2$, $4$ and $6$. That is, we have expressed a particular distribution as a 2-component mixture in two distinct ways, which means that $F(t) = \frac{1}{2}G_2(t-2) + \frac{1}{2}G_2(t-4)$ is not 2-identifiable. Yet, even without decomposing $F(t)$ in two distinct ways, we can immediately see that it cannot be 2-identifiable by noting that it is itself a symmetric distribution and therefore 1-identifiable (recall that no distribution can be $k$-identifiable for more than one $k$). Thus, only asymmetric elements of $\mathcal{M}_2$ can be 2-identifiable. We prove in Theorem 2 that asymmetry is not only necessary but also sufficient.

As expressed by equation (3), a location mixture may be written as a convolution. By exploiting the fact that convolution corresponds to multiplication of characteristic functions, we prove in Appendix A the following simple characterization of the set $\Omega_k^*$. Recall that $\Delta_k^-(\boldsymbol{\lambda}, \boldsymbol{\mu})$ denotes the reflection of $\Delta_k(\boldsymbol{\lambda}, \boldsymbol{\mu})$ about the origin.

THEOREM 1. *For $k \geq 1$, $\Omega_k^*$ defined in equation (4) is the set of $(\boldsymbol{\lambda}, \boldsymbol{\mu}) \in \Omega_k$ such that $\Delta_k^-(\boldsymbol{\lambda}, \boldsymbol{\mu})$ is the unique $k$-point distribution that yields a zero-symmetric distribution when convolved with $\Delta_k(\boldsymbol{\lambda}, \boldsymbol{\mu})$.*

It remains to describe $\Omega_k^*$ explicitly for certain values of $k$. It is not difficult to see that $\Omega_1^* = \Omega_1$. For the case $k = 2$, $\Omega_2^*$ cannot contain any $(\boldsymbol{\lambda}, \boldsymbol{\mu})$ for which $\lambda_1 = 0$, $\lambda_1 = \frac{1}{2}$ or $\lambda_1 = 1$ since those values make the mixture itself symmetric. But the symmetric mixtures in the case $k = 2$ are the only ones that are not 2-identifiable, as Theorem 2 states.

THEOREM 2. *$\Omega_2^* = \{(\boldsymbol{\lambda}, \boldsymbol{\mu}) \in \Omega_2 : \lambda_1 \notin \{0, \frac{1}{2}, 1\}\}$. Furthermore, every 2-identifiable mixture $F \in \mathcal{M}_2$ can be expressed as $G \star \Delta_2(\boldsymbol{\lambda}, \boldsymbol{\mu})$ for some $(\boldsymbol{\lambda}, \boldsymbol{\mu}) \in \Omega_2^*$ and $G \in \mathcal{S}$.*

The second statement in Theorem 2 may appear trivial at first glance. Yet, it is not immediately clear for general $k$ whether there exist $(\boldsymbol{\lambda}', \boldsymbol{\mu}') \in \Omega_k \setminus \Omega_k^*$ and $G \in \mathcal{S}$ such that $G \star \Delta_k(\boldsymbol{\lambda}', \boldsymbol{\mu}')$ is $k$-identifiable. Given an arbitrary $G \in \mathcal{S}$, the definition of $\Omega_k^*$ merely states that $(\boldsymbol{\lambda}, \boldsymbol{\mu}) \in \Omega_k^*$ is a sufficient condition, not a necessary condition, that $G \star \Delta_k(\boldsymbol{\lambda}, \boldsymbol{\mu})$ be $k$-identifiable. Theorem 2 states that for $k = 2$, it is also a necessary condition.

Unfortunately, the situation is not as straightforward when $k > 2$ as it is when $k = 2$. For instance, the set $\Omega_k^*$ does not contain all $(\boldsymbol{\lambda}, \boldsymbol{\mu}) \in \Omega_k$ such that $\Delta_k(\boldsymbol{\lambda}, \boldsymbol{\mu})$ is asymmetric and all $\lambda_j$ are nonzero. Furthermore, we do not know whether there exist $k$-identifiable distributions $F$ such that $\varphi_k^{-1}(F)$ is not in $\Omega_k^* \times \mathcal{S}$. Because it is somewhat complicated, we have put the statement of the explicit form of $\Omega_3^*$ into Appendix A as Theorem A.1. Here, we offer only the following sufficient (but certainly not necessary) condition for membership in $\Omega_3^*$.

COROLLARY 1. $(\boldsymbol{\lambda}, \boldsymbol{\mu}) \in \Omega_3^*$ if $\lambda_1 \lambda_2 \lambda_3 \neq 0$ and $(\mu_2 - \mu_1)/(\mu_3 - \mu_2) \notin \{\frac{1}{3}, \frac{1}{2}, 1, 2, 3\}$.

The stipulations in Corollary 1 that $\lambda_1 \lambda_2 \lambda_3 \neq 0$ and $(\mu_2 - \mu_1)/(\mu_3 - \mu_2) \neq 1$ ensure that $\Delta_3(\boldsymbol{\lambda}, \boldsymbol{\mu})$ cannot itself be symmetric; the stipulation that the larger of $\mu_2 - \mu_1$ and $\mu_3 - \mu_2$ cannot be two or three times the smaller eliminates two troublesome special cases (the only two, it turns out) which are given in equations (A.2)–(A.5).

Theorem 2 and Corollary 1 together imply that for $k \leq 3$, a $k$-component mixture of location-shifted symmetric distributions is almost always $k$-identifiable, in the sense that the set $\Omega_k \setminus \Omega_k^*$ has Lebesgue measure zero in $\mathbb{R}^{2k-1}$ (because of the constraint on $\boldsymbol{\lambda}$, $\Omega_k$ should be viewed as a subset of $\mathbb{R}^{2k-1}$, rather than $\mathbb{R}^{2k}$, in order to have positive Lebesgue measure). We conjecture that this is true for all $k$; however, because the situation gets even more complicated for larger $k$, we do not describe $\Omega_k^*$ for $k \geq 4$ in this article.

**3. Estimation.** Given a simple random sample from the distribution $F_0 = G_0 \star \Delta_k(\boldsymbol{\lambda}^0, \boldsymbol{\mu}^0)$, assuming $(\boldsymbol{\lambda}^0, \boldsymbol{\mu}^0)$ is contained in $\Omega_k^*$, it is natural to ask how one might tackle the problem of deconvolution. The idea for estimating $(\boldsymbol{\lambda}^0, \boldsymbol{\mu}^0)$ is as follows. Since a distribution is zero-symmetric if and only if its convolution with $G_0$ is zero-symmetric, Theorem 1 implies that there is exactly one $(\boldsymbol{\lambda}, \boldsymbol{\mu}) \in \Omega_k$ such that

$$F_0 \star \Delta_k^-(\boldsymbol{\lambda}, \boldsymbol{\mu}) = [G_0 \star \Delta_k(\boldsymbol{\lambda}^0, \boldsymbol{\mu}^0)] \star \Delta_k^-(\boldsymbol{\lambda}, \boldsymbol{\mu}) = G_0 \star [\Delta_k(\boldsymbol{\lambda}^0, \boldsymbol{\mu}^0) \star \Delta_k^-(\boldsymbol{\lambda}, \boldsymbol{\mu})]$$

is zero-symmetric, namely the true parameter value $(\boldsymbol{\lambda}^0, \boldsymbol{\mu}^0)$. Therefore, our plan is to search for a $\boldsymbol{\lambda}$ and a $\boldsymbol{\mu}$ that bring $\hat{F}_n \star \Delta_k^-(\boldsymbol{\lambda}, \boldsymbol{\mu})$ as close as possible to being a zero-symmetric distribution, where $\hat{F}_n$ is the empirical distribution function derived from the sample. We measure closeness to zero-symmetry as the distance between a distribution and its reflection across the origin. To this end, we define real-valued functions

$$(5) \qquad d(\boldsymbol{\lambda}, \boldsymbol{\mu}) = \mathcal{D}\{F_0 \star \Delta_k^-(\boldsymbol{\lambda}, \boldsymbol{\mu}), F_0^- \star \Delta_k(\boldsymbol{\lambda}, \boldsymbol{\mu})\}$$

and

$$(6) \qquad d_n(\boldsymbol{\lambda}, \boldsymbol{\mu}) = \mathcal{D}\{\hat{F}_n \star \Delta_k^-(\boldsymbol{\lambda}, \boldsymbol{\mu}), \hat{F}_n^- \star \Delta_k(\boldsymbol{\lambda}, \boldsymbol{\mu})\},$$

where $\mathcal{D}(F_1, F_2)$ is some measure of the distance between distribution functions $F_1$ and $F_2$. Provided $\mathcal{D}(F_1, F_2) = 0$ if and only if $F_1$ coincides with $F_2$, $(\boldsymbol{\lambda}^0, \boldsymbol{\mu}^0)$ is the unique minimizer of $d(\boldsymbol{\lambda}, \boldsymbol{\mu})$ and we estimate it by

$$(7) \qquad (\hat{\boldsymbol{\lambda}}, \hat{\boldsymbol{\mu}}) = \arg \min_{(\boldsymbol{\lambda}, \boldsymbol{\mu}) \in \Omega_k} d_n(\boldsymbol{\lambda}, \boldsymbol{\mu}).$$

[If the minimizer is not unique, then $(\hat{\boldsymbol{\lambda}}, \hat{\boldsymbol{\mu}})$ may be taken to be an arbitrarily selected minimizer.]

There are many possible choices for $\mathcal{D}$. Some, however, do not work well in this context. For example, total variation distance is not a good choice because $\hat{F}_n \star \Delta_k^-(\boldsymbol{\lambda}, \boldsymbol{\mu})$ and $\hat{F}_n^- \star \Delta_k(\boldsymbol{\lambda}, \boldsymbol{\mu})$ are both discrete distributions that are generally supported on entirely different points. In this article, we focus on $L_p$-distance for $1 \le p \le \infty$, defined for finite $p$ by

$$(8) \qquad \mathcal{D}(F_1, F_2) = \left( \int_{-\infty}^{\infty} |F_1(t) - F_2(t)|^p \, dt \right)^{1/p}$$

and for $p = \infty$ by $\mathcal{D}(F_1, F_2) = \sup_t |F_1(t) - F_2(t)|$. Note that if $p < \infty$, then $\mathcal{D}(F_1, F_2) < \infty$ whenever both $F_1$ and $F_2$ have finite first moments because

$$\{\mathcal{D}(F_1, F_2)\}^p \le \int_{-\infty}^{\infty} |F_1(t) - F_2(t)| \, dt$$

$$\le \int_{-\infty}^{0} \{F_1(t) + F_2(t)\} \, dt + \int_{0}^{\infty} \{1 - F_1(t) + 1 - F_2(t)\} \, dt$$

$$= \mathrm{E}_{F_1}|X| + \mathrm{E}_{F_2}|X|.$$

With this choice of distance and finite $p$,

$$(9) \qquad \{d(\boldsymbol{\lambda}, \boldsymbol{\mu})\}^p = \int_{-\infty}^{\infty} \left| \sum_{j=1}^{k} \lambda_j \{1 - F_0(\mu_j - t) - F_0(\mu_j + t)\} \right|^p dt$$

and

$$(10) \qquad \{d_n(\boldsymbol{\lambda}, \boldsymbol{\mu})\}^p = \int_{-\infty}^{\infty} \left| \sum_{j=1}^{k} \lambda_j \{1 - \hat{F}_n(\mu_j - t) - \hat{F}_n(\mu_j + t)\} \right|^p dt.$$

If $G_0(z)$ has finite first moment, then the minimizer $(\hat{\boldsymbol{\lambda}}, \hat{\boldsymbol{\mu}})$ of $d_n(\boldsymbol{\mu}, \boldsymbol{\lambda})$ is strongly consistent. To demonstrate this, we rely on a pair of lemmas.

LEMMA 1. *If $1 \le p < \infty$, we assume $G_0(z)$ has finite first moment; if $p = \infty$, we make no such assumption. In either case, $d(\hat{\boldsymbol{\lambda}}, \hat{\boldsymbol{\mu}}) \to 0$ almost surely as $n \to \infty$.*

LEMMA 2. *Under the assumptions of Lemma 1, for any $\varepsilon > 0$, there exists $\delta > 0$ such that $d(\boldsymbol{\lambda}, \boldsymbol{\mu}) > \delta$ whenever $\|(\boldsymbol{\lambda}, \boldsymbol{\mu}) - (\boldsymbol{\lambda}^0, \boldsymbol{\mu}^0)\| > \varepsilon$.*

Intuitively, Lemma 2 states that $d(\boldsymbol{\lambda}, \boldsymbol{\mu})$ is bounded away from zero outside any neighborhood of $(\boldsymbol{\lambda}^0, \boldsymbol{\mu}^0)$. By Lemma 2, the event $\{d(\hat{\boldsymbol{\lambda}}, \hat{\boldsymbol{\mu}}) < \delta\}$ is contained in the event $\{\|(\hat{\boldsymbol{\lambda}}, \hat{\boldsymbol{\mu}}) - (\boldsymbol{\lambda}^0, \boldsymbol{\mu}^0)\| \le \varepsilon\}$. But $\varepsilon$ is arbitrary, so by Lemma 1 we conclude that $\|(\hat{\boldsymbol{\lambda}}, \hat{\boldsymbol{\mu}}) - (\boldsymbol{\lambda}^0, \boldsymbol{\mu}^0)\| \to 0$ almost surely. This proves the following theorem.

THEOREM 3. *Suppose that $G_0(z)$ has finite first moment and $\mathcal{D}(\cdot, \cdot)$ is $L_p$-distance with $1 \le p \le \infty$. Then $(\hat{\boldsymbol{\lambda}}, \hat{\boldsymbol{\mu}}) \to (\boldsymbol{\lambda}^0, \boldsymbol{\mu}^0)$ almost surely as $n \to \infty$. (In the case $p = \infty$, the first moment condition is not necessary.)*

Once we have an estimate of $(\boldsymbol{\lambda}^0, \boldsymbol{\mu}^0)$, we turn to the question of estimating $G_0$. It may be that we only wish to estimate a particular functional of $G_0$ such as its variance $\sigma^2$. Since $F_0 = G_0 \star \Delta_k(\boldsymbol{\lambda}^0, \boldsymbol{\mu}^0)$, we obtain

$$\sigma^2 = \operatorname{Var}_{F_0}(X) - \sum_{j=1}^{k} \lambda_j (\mu_j - \bar{\mu})^2,$$

where $\bar{\mu} = \sum_j \lambda_j \mu_j$. In the case $k = 2$, if $S^2$ denotes the sample variance of $X_1, \dots, X_n$, then we obtain as an estimator of $\sigma^2$

$$(11) \qquad \hat{\sigma}^2 = S^2 - \hat{\lambda}_1 \hat{\lambda}_2 (\hat{\mu}_2 - \hat{\mu}_1)^2.$$

If $\sigma^2 < \infty$, then the strong consistency of $\hat{\sigma}^2$ follows from the strong law of large numbers and Theorem 3.

However, we may be interested in estimating the function $G_0(t)$ itself. We focus on the case $k = 2$. From $F_0 = G_0 \star \Delta_2(\boldsymbol{\lambda}^0, \boldsymbol{\mu}^0)$, we obtain the linear equation

$$(12) \qquad \begin{pmatrix} F_0(x) \\ F_0^-(x - \mu_1^0 - \mu_2^0) \end{pmatrix} = \begin{pmatrix} \lambda_1^0 & \lambda_2^0 \\ \lambda_2^0 & \lambda_1^0 \end{pmatrix} \begin{pmatrix} G_0(x - \mu_1^0) \\ G_0(x - \mu_2^0) \end{pmatrix},$$

valid for all $x$ for which $\mu_1^0 + \mu_2^0 - x$ is a continuity point of $F_0(\cdot)$. Equation (12) may easily be inverted to give a formula for $G_0(x - \mu_1^0)$ and $G_0(x - \mu_2^0)$ [note that the identifiability requirement that $\lambda_1 \ne 1/2$ is reflected in the fact that the $2 \times 2$ matrix in equation (12) is singular when $\lambda_1 = 1/2$]. Replacing the parameters $\boldsymbol{\lambda}^0$, $\boldsymbol{\mu}^0$ and $F_0$ by their respective estimates gives

$$\begin{pmatrix} G_0(x - \hat{\mu}_1) \\ G_0(x - \hat{\mu}_2) \end{pmatrix} \approx \frac{1}{\hat{\lambda}_1 - \hat{\lambda}_2} \begin{pmatrix} \hat{\lambda}_1 \hat{F}_n(x) - \hat{\lambda}_2 \hat{F}_n^-(x - \hat{\mu}_1 - \hat{\mu}_2) \\ -\hat{\lambda}_2 \hat{F}_n(x) + \hat{\lambda}_1 \hat{F}_n^-(x - \hat{\mu}_1 - \hat{\mu}_2) \end{pmatrix}.$$

We could thus obtain two estimates for $G_0(z)$, one by setting $x - \hat{\mu}_1 = z$ and the other by setting $x - \hat{\mu}_2 = z$. Taking the mean of these two estimates yields

$$(13) \quad \hat{G}_0(z) = \frac{\hat{\lambda}_1 [\hat{F}_n^-(z - \hat{\mu}_1) + \hat{F}_n(z + \hat{\mu}_1)] - \hat{\lambda}_2 [\hat{F}_n^-(z - \hat{\mu}_2) + \hat{F}_n(z + \hat{\mu}_2)]}{2(\hat{\lambda}_1 - \hat{\lambda}_2)}.$$

The function $\hat{G}_0(z)$ has the appealing property that it satisfies the zero-symmetry condition: at all continuity points $z$, $\hat{G}_0(z) = 1 - \hat{G}_0(-z)$. Furthermore, $\lim_{z \to -\infty} \hat{G}_0(z) = 0$ and $\lim_{z \to \infty} \hat{G}_0(z) = 1$. However, $\hat{G}_0(z)$ is not necessarily a legitimate distribution function because it is not generally nondecreasing. Although this may initially appear to be a drawback, we consider the following corollary of Theorem 3 and the Glivenko–Cantelli theorem ([2], page 275), which states that $\sup_t |\hat{F}_n(t) - F_0(t)| \to 0$ almost surely.

COROLLARY 2. *Under the assumptions of Theorem 3 with $k = 2$, $\sup_z |\hat{G}_0(z) - G_0(z)| \to 0$ almost surely as $n \to \infty$.*

Thus, if we compute and graph an estimate $\hat{G}_0(z)$ and find that it is not roughly monotone increasing, then there are two possible causes: either the sample size is too small for the asymptotics of Corollary 2 or the model is misspecified. In other words, $\hat{G}_0(z)$ might serve as a sort of graphical goodness-of-fit test; we will say more about this in Section 5. For $k > 2$, derivation of an estimator of $G_0(x)$ is not as straightforward as for $k = 2$ since $\hat{G}_0$ may not be easily attained as the solution of a system of linear equations; a different method of deconvolution may be necessary in this case.

**4. Generalizing the Hodges–Lehmann estimator.** Although the strong consistency proved in Theorem 3 and the resulting Corollary 2 are valid for $L_p$-distance for any $1 \le p \le \infty$, this section and the next consider only $p = 2$. Here, we demonstrate that the proposed estimator of equation (7), where $\mathcal{D}$ is $L_2$-distance as defined in equation (8), is a generalization of the Hodges–Lehmann estimator (1) to the case of finite mixtures. Furthermore, we establish sufficient conditions for the asymptotic normality of the estimator when $p = 2$.

Let $H_W(t)$ denote the distribution function of an arbitrary random variable $W$. Suppose $W$, $W_1$ and $W_2$ are independent and identically distributed random variables. Then denoting $\max\{W_1, W_2\}$ by $W_1 \vee W_2$, the identities $H_W^2(t) = H_{W_1 \vee W_2}(t)$ and $H_W(t)H_{-W}(t) = H_{-W_1 \vee W_2}(t)$ imply that

$$\int_{-\infty}^{\infty} \{H_W(t) - H_{-W}(t)\}^2 \, dt$$

$$= \int_{-\infty}^{0} \{H_{W_1 \vee W_2}(t) + H_{-W_1 \vee -W_2}(t) - 2H_{-W_1 \vee W_2}(t)\} \, dt$$

$$- \int_{0}^{\infty} \{-H_{W_1 \vee W_2}^-(-t) - H_{-W_1 \vee -W_2}^-(-t) + 2H_{-W_1 \vee W_2}^-(-t)\} \, dt$$

$$= \mathrm{E}\{2(-W_1 \vee W_2) - (W_1 \vee W_2) - (-W_1 \vee -W_2)\}.$$

Since $2\max\{a, b\} = a + b + |a - b|$, we obtain

$$(14) \qquad \int_{-\infty}^{\infty} \{H_W(t) - H_{-W}(t)\}^2 \, dt = \mathrm{E}(|W_1 + W_2| - |W_1 - W_2|).$$

Letting $W \sim \hat{F}_n \star \Delta_k^-(\boldsymbol{\lambda}, \boldsymbol{\mu})$, we may combine equation (6) with equation (14) to obtain

$$(15) \qquad \begin{aligned} \{d_n(\boldsymbol{\theta})\}^2 &= \mathrm{E}(|W_1 + W_2| - |W_1 - W_2|) \\ &= \frac{1}{n^2} \sum_{i=1}^{n} \sum_{j=1}^{n} f_{\boldsymbol{\theta}}(x_i, x_j), \end{aligned}$$

where $\boldsymbol{\theta} = (\boldsymbol{\lambda}, \boldsymbol{\mu})$ and

$$(16) \qquad f_{\boldsymbol{\theta}}(x_i, x_j) = \sum_{a=1}^{k} \sum_{b=1}^{k} \lambda_a \lambda_b (|x_i + x_j - \mu_a - \mu_b| + |x_i - x_j - \mu_a + \mu_b|).$$

When $k = 1$, the only parameter to estimate is $\mu$, the center of symmetry, and minimization of expression (15) reduces to

$$(17) \qquad \hat{\mu} = \arg\min_{\mu} \sum_{i=1}^{n} \sum_{j=1}^{n} \left| \frac{x_i + x_j}{2} - \mu \right|.$$

Comparing $\hat{\mu}$ with the Hodges–Lehmann estimator $\hat{\mu}_{\mathrm{HL}}$ of equation (1), the two estimators are nearly the same, except that the sum in (1) places twice as much weight on the cases when $i = j$. Based on this similarity when $k = 1$, the estimator

$$(18) \qquad \hat{\boldsymbol{\theta}} = \arg\min_{\boldsymbol{\theta}} d_n(\boldsymbol{\theta}) = \arg\min_{\boldsymbol{\theta}} \frac{1}{n^2} \sum_{i=1}^{k} \sum_{j=1}^{k} f_{\boldsymbol{\theta}}(x_i, x_j)$$

may be categorized as a generalization of the Hodges–Lehmann estimator for $k \geq 2$.

To establish the asymptotic normality of $\hat{\boldsymbol{\theta}}$, note that $\{d_n(\boldsymbol{\theta})\}^2$ from equation (15) is a $V$-process (i.e., a set of $V$-statistics indexed by the parameter $\boldsymbol{\theta}$). Define the functionals

$$V^{(k)}(h) = \mathrm{E}h(X_1, \ldots, X_k)$$

and

$$V_n^{(k)}(h) = \frac{1}{n^k} \sum_{i_1=1}^{n} \cdots \sum_{i_k=1}^{n} h(X_{i_1}, \ldots, X_{i_k}),$$

where $h$ is some scalar- or vector-valued function of $k$ real variables (recall that $X_1, \ldots, X_n$ is a random sample from $F_0$). In particular, $\{d(\boldsymbol{\theta})\}^2 = V^{(2)}(f_{\boldsymbol{\theta}})$ and $\{d_n(\boldsymbol{\theta})\}^2 = V_n^{(2)}(f_{\boldsymbol{\theta}})$. The Hoeffding decomposition for $V$-statistics has exactly the same form as it has for $U$-statistics (cf. [1]),

$$(19) \qquad V_n^{(m)}(h) - V^{(m)}(h) = \sum_{k=1}^{m} \binom{m}{k} V_n^{(k)}(\pi_k h),$$

where $\pi_k h$ is the $k$th Hoeffding projection defined, using the notation of empirical processes, by

$$(20) \qquad \pi_k h(x_1, \ldots, x_k) = (\delta_{x_1} - F_0) \cdots (\delta_{x_k} - F_0) F_0^{m-k} h,$$

where $Qf = \int f \, dQ$ denotes the action of the expectation operator under the distribution $Q$ on the function $f$ and $\delta_x$ denotes a point mass at $x$ (see [9] for an alternative formulation of the $\pi_k h$ projection functions). Therefore, the sufficient conditions established by Arcones, Chen and Giné [1] for the asymptotic normality of $U$-processes are valid in the present case. Their Theorem 2.1, establishing asymptotic normality, is based on Theorem 2 of [12] and we adapt these theorems to the present situation as follows.

THEOREM 4. *Assume that $G_0(z)$ has finite first moment and that the following hold:*

(i) *$V^{(2)}(f_{\boldsymbol{\theta}})$, as a function of $\boldsymbol{\theta}$, has strictly positive definite second derivative $J$ at its minimizing value $\boldsymbol{\theta}^0$;*

(ii) *for any $\varepsilon > 0$, there exists $\delta > 0$ such that*

$$\limsup_{n \to \infty} P\left\{ \sup_{\boldsymbol{\theta} \in B_\delta} |nV_n^{(2)}(\pi_2 f_{\boldsymbol{\theta}} - \pi_2 f_{\boldsymbol{\theta}^0})| > \varepsilon \right\} < \varepsilon,$$

*where $B_\delta$ denotes the open ball of radius $\delta$ centered at $\boldsymbol{\theta}^0$;*

(iii) *there exists a measurable vector-valued function $\Delta(x)$ satisfying $\mathrm{E}\Delta(X) = \mathbf{0}$, $\mathrm{E}\|\Delta(X)\|^2 < \infty$ and*

$$\pi_1 f_{\boldsymbol{\theta}}(x) = \pi_1 f_{\boldsymbol{\theta}^0}(x) + (\boldsymbol{\theta} - \boldsymbol{\theta}^0)^t \Delta(x) + \|\boldsymbol{\theta} - \boldsymbol{\theta}^0\| r_{\boldsymbol{\theta}}(x)$$

*for some $r_{\boldsymbol{\theta}}$ such that for any $\varepsilon > 0$, there exists $\delta > 0$ such that*

$$\limsup_{n \to \infty} P\Big\{ \sup_{\boldsymbol{\theta} \in B_\delta} |\sqrt{n} V_n^{(1)}(r_{\boldsymbol{\theta}})| > \varepsilon \Big\} < \varepsilon,$$

*where $B_\delta$ denotes the open ball of radius $\delta$ centered at $\boldsymbol{\theta}^0$.*

*Then $\sqrt{n}(\hat{\boldsymbol{\theta}} - \boldsymbol{\theta}^0) \to N(0, 4J^{-1}\Sigma J^{-1})$ in distribution, where $\Sigma$ is the covariance matrix of $\Delta(X)$.*

We caution that although the covariance formula concluding Theorem 4 appears simple, in our experience, it is extremely complicated to use in practice. Thus, we recommend a bootstrapped estimate of the estimator's covariance if one is needed.

The general idea of the proof of Theorem 4 is as follows. First, define $\tilde{\boldsymbol{\theta}} = \boldsymbol{\theta}^0 - 2J^{-1}V_n^{(1)}\Delta$. Show that both $\sqrt{n}(\hat{\boldsymbol{\theta}} - \boldsymbol{\theta}^0)$ and $\sqrt{n}(\tilde{\boldsymbol{\theta}} - \boldsymbol{\theta}^0)$ are stochastically bounded (the latter fact follows directly from the central limit theorem). Use these facts to prove that

$$(21) \qquad nV_n^{(2)}(f_{\hat{\boldsymbol{\theta}}} - f_{\boldsymbol{\theta}^0}) = -\frac{n}{2}(\tilde{\boldsymbol{\theta}} - \boldsymbol{\theta}^0)^t J(\tilde{\boldsymbol{\theta}} - \boldsymbol{\theta}^0) + o_P(1),$$

$$(22) \qquad \begin{aligned} nV_n^{(2)}(f_{\hat{\boldsymbol{\theta}}} - f_{\boldsymbol{\theta}^0}) &= \frac{n}{2}(\hat{\boldsymbol{\theta}} - \boldsymbol{\theta}^0)^t J(\hat{\boldsymbol{\theta}} - \boldsymbol{\theta}^0) \\ &\quad - n(\hat{\boldsymbol{\theta}} - \boldsymbol{\theta}^0)^t J(\tilde{\boldsymbol{\theta}} - \boldsymbol{\theta}^0) + o_P(1). \end{aligned}$$

Now, since $\hat{\boldsymbol{\theta}}$ minimizes $V_n^{(2)}f_{\boldsymbol{\theta}}$, subtract equation (22) from equation (21) to obtain

$$0 \le -\frac{n}{2}\|J^{1/2}(\hat{\boldsymbol{\theta}} - \tilde{\boldsymbol{\theta}})\| + o_P(1),$$

which implies that $\sqrt{n}(\hat{\boldsymbol{\theta}} - \tilde{\boldsymbol{\theta}}) \to 0$ in probability. Since $\sqrt{n}(\tilde{\boldsymbol{\theta}} - \boldsymbol{\theta}^0)$ converges in distribution to $N(0, 4J^{-1}\Sigma J^{-1})$ by the central limit theorem, this proves the result. Since the technical details missing above do not differ from those in the proof of Theorem 2.1 in [1], we do not repeat them here. These arguments are based on the proofs of Theorem 2 and Lemma 3 in [12], to which we also refer the interested reader.

## 5. Implementation and examples.

Although equation (15) allows an intuitive appreciation of the estimation method based on $L_2$-norm minimization, it is not the most convenient formula from a computational standpoint. To aid notation, we introduce the functional inner product

$$\langle f, g \rangle = \int_{-\infty}^{\infty} f(t)g(t)\, dt$$

and let $\|f\| = \sqrt{\langle f, f \rangle}$ denote the corresponding norm. Therefore, if we define

$$a_j(t; \boldsymbol{\mu}) = \frac{1}{n} \sum_{i=1}^{n} I\{\mu_j - x_i \le t\} - \frac{1}{n} \sum_{i=1}^{n} I\{x_i - \mu_j \le t\}$$

for $j = 1, \ldots, k$, then equation (10) implies that we may succinctly write $d_n(\boldsymbol{\lambda}, \boldsymbol{\mu}) = \| \sum_j \lambda_j a_j \|$.

In the case $k = 2$, we have $\lambda_1 + \lambda_2 = 1$ and thus $d_n(\boldsymbol{\lambda}, \boldsymbol{\mu}) = \|\lambda_1 a_1 + \lambda_2 a_2\|^2$ is a quadratic function in $\lambda_1$, minimized by $\hat{\lambda}_1 = -\langle a_1 - a_2, a_2 \rangle / \|a_1 - a_2\|^2$. Substituting $\hat{\lambda}_1$ into $d_n(\boldsymbol{\lambda}, \boldsymbol{\mu})$ gives

$$m(\mu_1, \mu_2) = \frac{\|a_1\|^2 \|a_2\|^2 - \langle a_1, a_2 \rangle^2}{\|a_1 - a_2\|^2}$$

as the function we wish to minimize. We accomplish this minimization using the "optim" function in R, which implements an iterative Nelder–Mead algorithm [type "help(optim)" in R for details]. We recommend multiple starting values due to the fact that $m(\boldsymbol{\mu})$ often has multiple local minima. Code that implements our method for $k = 2$, written in R, is available at www.stat.psu.edu/˜dhunter/code.

We first apply our semi-parametric estimation procedure (hereafter referred to as SP) to the Old Faithful data depicted in Figure 1. Results are compared with those obtained by maximizing the two-component normal mixture likelihood (a procedure we call NMLE from now on) that assumes components with equal variances. The assumption of equal variances in the NMLE case not only provides a fair comparison with the SP method (which assumes components with exactly the same shape) but it also avoids the awkward situation of an unbounded likelihood function created when unequal variances are assumed in the normal mixture (we say more about this in Section 6).

In Table 1, we see very close agreement between the two methods, with the standard errors only slightly larger for the SP method, even though Figure 2 indicates that the data in each component appear to follow a normal distribution quite closely. This figure depicts the close agreement between the estimate $\hat{G}_0$ of equation (13) and the estimate based on the NMLE method, namely, a normal distribution function with mean 0 and variance 34.45. One important benefit of the SP method, even in cases like this in which the data appear to be well modeled by a normal distribution, is that we may validate the normality assumption using this nonparametric estimate of the underlying symmetric cdf.

Next, we compare the SP method to the NMLE method for various types of simulated datasets. Results are summarized in Figure 3. The symmetric component distributions are taken to be normal, double exponential, uniform or $t_2$ ($t$ on two degrees of freedom). The values of $\lambda_1^0$ are taken

TABLE 1

*Shown here are parameter estimates along with bootstrapped standard errors for the semi-parametric approach (SP) and the normal mixture approach using maximum likelihood estimation (NMLE) for the Old Faithful geyser data. The bootstrapped estimates are based on 200 resamples*

|      | $\hat{\mu}_1$ (SE) | $\hat{\mu}_2$ (SE) | $\hat{\lambda}$ (SE) | $\hat{\sigma}^2$ (SE) |
|------|--------------------|--------------------|----------------------|-----------------------|
| SP   | 54.00 (0.76)       | 80.00 (0.50)       | 0.352 (0.032)        | 30.66 (7.93)          |
| NMLE | 54.61 (0.67)       | 80.09 (0.45)       | 0.361 (0.032)        | 34.45 (3.39)          |

from the set $\{0.15, 0.30, 0.45\}$ and each sample is of size $n = 200$. Because both methods are susceptible to finding local optimum points, each algorithm was started at several places for every dataset. The initial values of $\mu_1$ and $\mu_2$ were taken to be the $q_1$ and $q_2$ sample quantiles for each of the ten possible combinations satisfying $q_1, q_2 \in \{0.05, 0.2, 0.5, 0.8, 0.95\}$ and $q_1 < q_2$. Furthermore, the EM algorithm for the NMLE method was started with initial values $\lambda_1 = 0.5$ and $\sigma^2$ equal to half of the sample variance. The parameter estimates in all cases were taken to be the values corresponding to the best value of the objective function (i.e., the lowest value of $d_n$ or the highest value of the likelihood) among the ten. We took $\mu_1^0 = -\mu_2^0 = -1$ and $\sigma_0^2 = 1$ for all normal, double exponential and uniform examples. For the $t_2$ distribution, which has infinite variance, we took $\mu_1^0 = -\mu_2^0 = -2$.

For heavy-tailed distributions such as the double exponential and especially $t_2$, the SP method outperforms the NMLE method. Perhaps, since $\lambda_1^0 = 0.30$ is the farthest value from the nonidentifiable situations $\lambda_1^0 = 0$ and $\lambda_1^0 = 0.5$, the SP method fares better at that value of $\lambda_1^0$ than at $\lambda_1^0 = 0.15$ or $\lambda_1^0 = 0.45$. Nonetheless, the SP method performed surprisingly well for the $\lambda_1^0 = 0.45$ case, despite the fact that 0.45 is so close to the nonidentifiable value of 0.5. Both SP and NMLE had the most difficult time at $\lambda_1^0 = 0.15$, regardless of the type of component distributions.

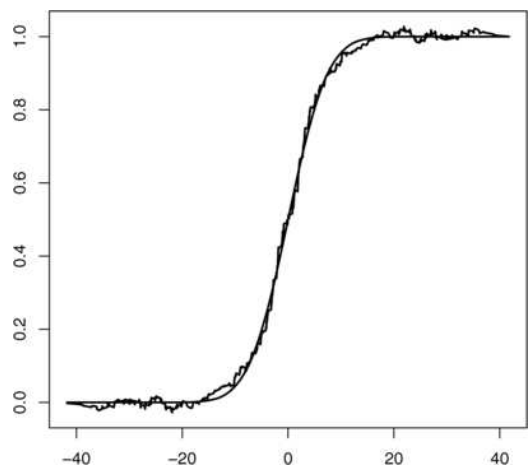

FIG. 2. *The jagged line is $\hat{G}_0$, estimated from the Old Faithful data using equation* (13), *and the other line is the NMLE estimate of $G_0$, which is forced to be normal.*

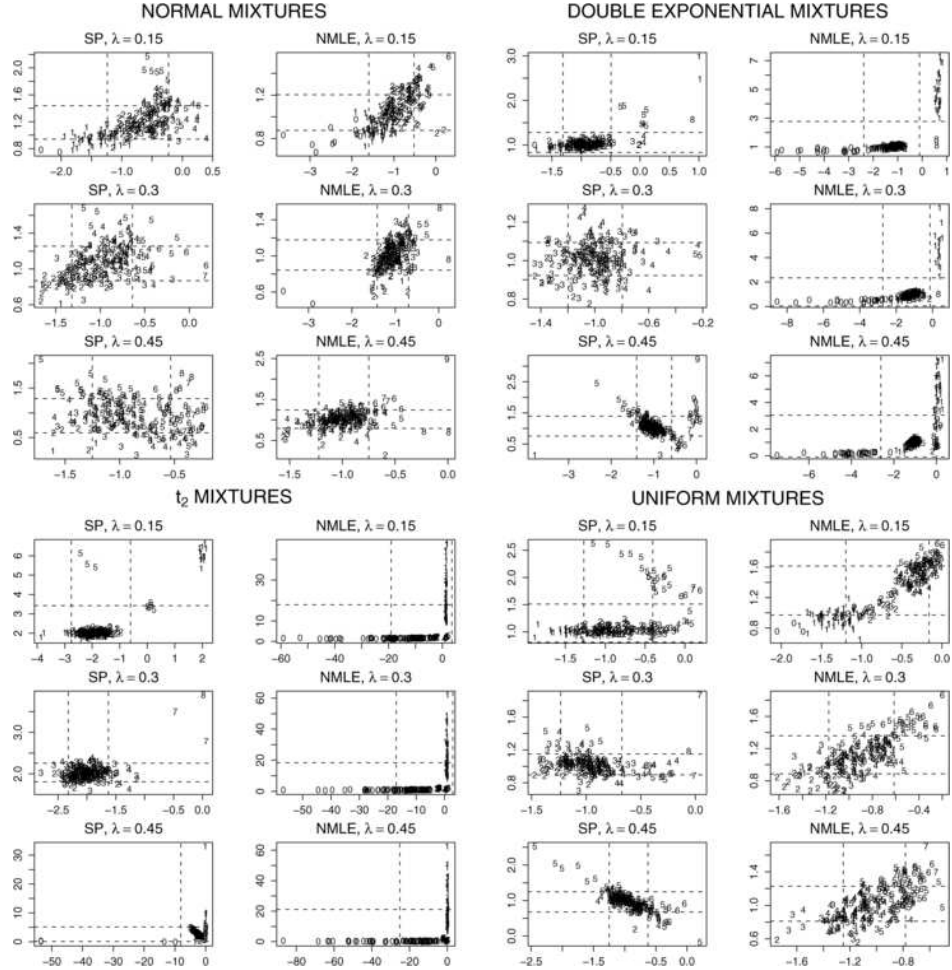

FIG. 3. *Shown above are scatterplots of the parameter estimates* $(\hat{\mu}_1, \hat{\mu}_2)$ *from 200 simulated datasets of size 200. The true value of* $(\mu_1^0, \mu_2^0)$ *is* $(-1, 1)$ *in all plots except the* $t_2$ *plots, where it is* $(-2, 2)$. *Each point is represented by the leading digit of* $\lambda$ *rounded to the nearest tenth. The dashed lines are one sample standard deviation of* $\hat{\mu}_i$ *on either side of the sample mean of* $\hat{\mu}_i$ *for* $i = 1, 2$.

Finally, we consider the robustness of both methods to violations of the assumptions that the component distributions are symmetric and that the components differ only in location. For nonsymmetric distributions, we use

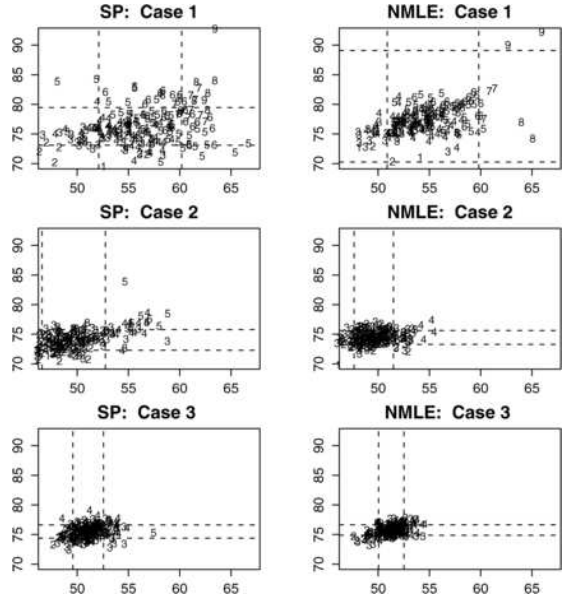

Fig. 4. *Estimation methods applied to simulated datasets that violate assumptions. All sets of axes are identical and the content of each plot is as explained in Figure 3. In case 1, both symmetry and equal variance assumptions are violated; in case 2, only symmetry is violated; in case 3, only equal variances is violated. In each case, $df_1 = 50$ and $df_2 = 75$.*

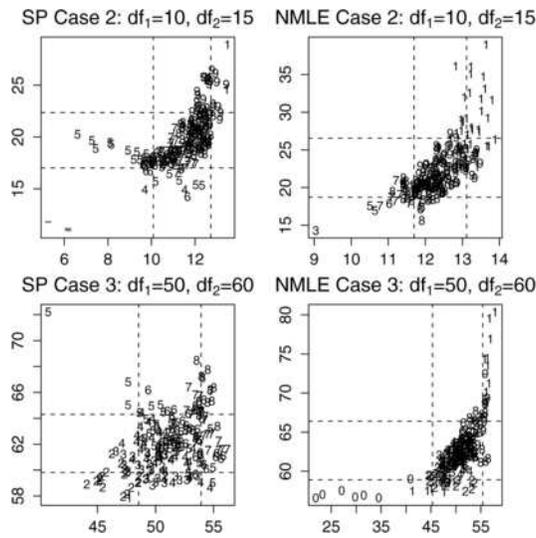

Fig. 5. *Both the upper plots and lower plots compare unfavorably with the values $df_1 = 50$ and $df_2 = 75$ used in Figure 4.*

$\chi^2$ distributions. Suppose that $df_1$ and $df_2$ are given positive integers, where $df_1 < df_2$. Then we consider the following three cases:

1. $\lambda\chi^2(df_1) + (1 - \lambda)\chi^2(df_2)$;
2. $\lambda\chi^2(df_1) + (1 - \lambda)[(df_2 - df_1) + \chi^2(df_1)]$;
3. $\lambda N(df_1, 2df_1) + (1 - \lambda)N(df_2, 2df_2)$.

Both assumptions are violated in case 1, only the symmetry assumption is violated in case 2 and only the equal-variances assumption is violated in case 3.

Consider $df_1 = 50$ and $df_2 = 75$ (all three cases) as an example (Figure 4). Not surprisingly, case 1 fares the worst. In a comparison between case 2 and case 3, it appears that case 3 is slightly better overall, suggesting that the skewness causes greater problems than the unequal variances. This impression is further strengthened by the fact that $df_1 = 10$, $df_2 = 15$ fares much worse in case 2 than $df_1 = 50$, $df_2 = 75$ (see the upper plots in Figure 5). On the other hand, poorly-separated means appears to be more detrimental than mismatched variances, since case 3 fares much better for $df_1$ and $df_2$ well separated than close together. For instance, with $df_1 = 50$, the performance improves in case 2 as $df_2$ increases, say, from 60 to 75 (see the lower plots in Figure 5). Further numerical tests reinforce these general impressions.

In summary, we find that for both methods, the presence of unequal variances appears to be the least detrimental violation of the assumptions; it is not even as serious as poorly-separated means. On the other hand, violation of the symmetry assumption appears to have a much stronger and more unpredictable effect on the results.

**6. Discussion.** This article addresses the question of how restrictive the assumptions about a mixture distribution must be in order for identifiability to hold and hence for estimation to make sense. In particular, we investigate the effect of presuming that the component distributions are symmetric and that they are all the same, except for location shifts. We establish comprehensive identifiability results in the 2- and 3-component cases and indicate a direction of analysis for $k > 3$ (although similar results appear to become prohibitively complicated for larger $k$). We emphasize that with so little assumed about the form of the component distributions, it is quite surprising that identifiability is provable at all. Relaxing our assumptions even further by, say, allowing location-scale transformations (instead of just location transformations) is likely to have a major impact on identifiability by undermining the convolution structure so central to the theoretical development in this article. Another possible assumption one might consider is unimodality, instead of symmetry, of components. However, such an assumption would certainly destroy identifiability unless additional restrictions were assumed;

furthermore, Walther [14, 15] has shown that mode-hunting and the search for mixture structure are not always compatible.

There are clear practical applications of the estimation method we propose. The simulation studies and data analysis in Section 5 suggest that our method is robust to different component distribution shapes and that its performance is never much worse than, and sometimes much better than, the canonical method of maximum likelihood assuming normal components. In addition, we have shown in a two-component example how our method allows the validation of parametric assumptions used to fit mixture models by providing a nonparametric estimate of the (unshifted) zero-symmetric components.

The assumption that each of the component distributions must have exactly the same shape, which means, among other things, that they must have the same variance if the variance is defined, may, at first glance, seem like an overly rigid assumption. In particular, when comparing our method with the method of maximum likelihood using normal components, it may seem that the normal approach has an advantage since it appears to offer the possibility of allowing different component variances. However, the equal-variance assumption is quite common throughout statistics (e.g., in ANOVA) and for cases in which this assumption is appropriate, it is wise—both in terms of statistical power and model parsimony—to utilize inference techniques that implement it. Furthermore, the comparison with the normal approach is slightly misleading: if we assume unknown means and unknown variances that are different, then the normal mixture likelihood is unbounded and therefore has no maximizer (although this defect may be mended by placing a positive lower bound on the component variances). Finally, we point out that our method is consistent, even when there is no finite second moment such as in the simulated $t$-distribution examples of Section 5, a case in which the parametric method performs poorly. We believe that the SP method offers an attractive alternative and/or complement to the standard parametric estimation in many problems.

## APPENDIX A: IDENTIFIABILITY PROOFS

We first prove Theorem 1. For general $k$, define independent random variables $Y$ and $Z$ such that $Z \sim G$ and $Y \sim \Delta_k(\boldsymbol{\lambda}, \boldsymbol{\mu})$, where $G \in \mathcal{S}$ and $(\boldsymbol{\lambda}, \boldsymbol{\mu}) \in \Omega_k$. Then by equation (3), $X = Y + Z$ has the mixture distribution $F$. In terms of the characteristic functions of these random variables, we have $\phi_X(t) = \phi_Y(t)\phi_Z(t)$. Suppose that

$$(\text{A.1}) \qquad \phi_Y(t)\phi_Z(t) = \phi_{Y'}(t)\phi_{Z'}(t)$$

for independent $Z' \sim G'$ and $Y' \sim \Delta_k(\boldsymbol{\lambda}', \boldsymbol{\mu}')$. Note that $k$-identifiability of $F$ holds if and only if equation (A.1) implies $Y \overset{\mathcal{D}}{=} Y'$ and $Z \overset{\mathcal{D}}{=} Z'$, where $\overset{\mathcal{D}}{=}$ denotes equality of distributions. A pair of lemmas will help us to prove the theorem.

LEMMA A.1.   *Equation* (A.1) *implies that* $Y - Y'$ *is zero-symmetric.*

Note that $Y - Y'$ can always be made zero-symmetric by taking $Y' \overset{\mathcal{D}}{=} Y$. If this choice of $Y'$ is the only choice that makes $Y - Y'$ zero-symmetric, then $k$-identifiability follows, as stated in Lemma A.2.

LEMMA A.2.   *For* $Y \sim \Delta_k(\boldsymbol{\lambda}, \boldsymbol{\mu})$, *suppose that* $Y - Y'$ *is zero-symmetric for independent* $Y' \sim \Delta_k(\boldsymbol{\lambda}', \boldsymbol{\mu}')$ *only if* $\boldsymbol{\lambda} = \boldsymbol{\lambda}'$ *and* $\boldsymbol{\mu} = \boldsymbol{\mu}'$. *Then for any* $Z \sim G \in \mathcal{S}$, *the distribution of* $X = Y + Z$ *is* $k$-*identifiable.*

Lemma A.2 implies that $\Omega_k^*$ must contain all $(\boldsymbol{\lambda}, \boldsymbol{\mu})$ such that $\Delta_k(\boldsymbol{\lambda}, \boldsymbol{\mu})$ cannot be zero-symmetrized by convolution with any $k$-point distribution other than $\Delta_k^-(\boldsymbol{\lambda}, \boldsymbol{\mu})$. But these $(\boldsymbol{\lambda}, \boldsymbol{\mu})$ are the only possible elements of $\Omega_k^*$, for if $\Delta_k(\boldsymbol{\lambda}, \boldsymbol{\mu}) \star \Delta_k^-(\boldsymbol{\lambda}', \boldsymbol{\mu}')$ is zero-symmetric, but $(\boldsymbol{\lambda}', \boldsymbol{\mu}') \neq (\boldsymbol{\lambda}, \boldsymbol{\mu})$, then

$$\Delta_k^-(\boldsymbol{\lambda}, \boldsymbol{\mu}) \star \{\Delta_k(\boldsymbol{\lambda}, \boldsymbol{\mu}) \star \Delta_k^-(\boldsymbol{\lambda}', \boldsymbol{\mu}')\} = \Delta_k^-(\boldsymbol{\lambda}', \boldsymbol{\mu}') \star \{\Delta_k(\boldsymbol{\lambda}, \boldsymbol{\mu}) \star \Delta_k^-(\boldsymbol{\lambda}, \boldsymbol{\mu})\}$$

is not $k$-identifiable. This proves Theorem 1—all that remains is to prove Lemmas A.1 and A.2.

PROOF OF LEMMA A.1.   A random variable is zero-symmetric if and only if its characteristic function is real-valued ([2], page 363). Multiplying each side of equation (A.1) by the complex conjugate of $\phi_{Y'}(t)$, namely $\phi_{-Y'}(t)$, we conclude that $\phi_Y(t)\phi_{-Y'}(t)$ is real-valued for all $t$ such that $\phi_Z(t) \neq 0$. Since $\phi_Z(t)$ is nonzero in a neighborhood of $t = 0$, the analytic function

$$\text{Im}\{\phi_Y(t)\phi_{-Y'}(t)\} = \sum_{i=1}^{k}\sum_{j=1}^{k} \lambda_i \lambda_j' \sin t(\mu_i - \mu_j')$$

equals zero on an open interval and must thus be identically zero on the whole real line. We conclude that if equation (A.1) holds, then $\phi_Y(t)\phi_{-Y'}(t)$ is real-valued, so $Y - Y'$ is zero-symmetric.   $\square$

PROOF OF LEMMA A.2.   By Lemma A.1, equation (A.1) implies that $\phi_Y(t) = \phi_{Y'}(t)$, so $\phi_Z(t) = \phi_{Z'}(t)$ whenever $\phi_Y(t) \neq 0$. But $\phi_Y(t)$ is an analytic function that is not identically zero, so $\{t : \phi_Y(t) = 0\}$ is a discrete set. For continuous functions $\phi_Z(t) = \phi_{Z'}(t)$ to agree outside a discrete set, they must be identical. Therefore, equation (A.1) implies both $Y \overset{\mathcal{D}}{=} Y'$ and $Z \overset{\mathcal{D}}{=} Z'$, so the distribution of $Y + Z$ is $k$-identifiable.   $\square$

Next, we prove Theorem 2 and then state and prove a similar characterization of $\Omega_3^*$ in Theorem A.1.

PROOF OF THEOREM 2. For any $(\boldsymbol{\lambda}, \boldsymbol{\mu}) \in \Omega_2$ with $2\lambda_1 \in \{0, 1, 2\}$, $\Delta_2(\boldsymbol{\lambda}, \boldsymbol{\mu})$ is symmetric and thus $G \star \Delta_2(\boldsymbol{\lambda}, \boldsymbol{\mu})$ cannot be 2-identifiable. Conversely, take $Y \sim \Delta_2(\boldsymbol{\lambda}, \boldsymbol{\mu})$ with $2\lambda_1 \notin \{0, 1, 2\}$ and suppose that $Y' \sim \Delta_2(\boldsymbol{\lambda}', \boldsymbol{\mu}')$ is independent of $Y$ with the property that $Y - Y'$ is zero-symmetric.

The largest and smallest values assumed by $Y - Y'$ must be opposites and receive the same weight by zero-symmetry, which implies that $\mu_2 - \mu_1' = \mu_2' - \mu_1$ and $\lambda_1 \lambda_2' = \lambda_1' \lambda_2$. Thus, $\boldsymbol{\lambda} = \boldsymbol{\lambda}'$. If $\boldsymbol{\mu} \neq \boldsymbol{\mu}'$, then neither $\mu_1 - \mu_1'$ nor $\mu_2 - \mu_2'$ can be zero, so these points must receive the same weight, leading to $\lambda_1 = \lambda_2$, a contradiction. We conclude that $Y$ and $Y'$ have the same distribution. $\square$

We now consider the case $k = 3$. We start by giving two distinct cases in which a nonsymmetric 3-point distribution may be nontrivially symmetrized by convolution with another 3-point distribution [recall that $\Delta_3(\boldsymbol{\lambda}, \boldsymbol{\mu})$ may always be trivially symmetrized by convolution with $\Delta_3^-(\boldsymbol{\lambda}, \boldsymbol{\mu})$]. We then assert in Theorem A.1 that these are the only such 3-point distributions.

CASE 1. For any real numbers $c$, $d$ and $r$ such that $d > 0$ and $r > 1$, let

(A.2) $\qquad \boldsymbol{\mu} = (c, c + 4d, c + 6d) \quad \text{and} \quad \boldsymbol{\lambda} \propto (r^2, r^2 - 1, r),$

(A.3) $\qquad \boldsymbol{\mu}' = (c + d, c + 3d, c + 5d) \quad \text{and} \quad \boldsymbol{\lambda}' \propto (r, r + 1, 1).$

Then

$$\Delta_3(\boldsymbol{\lambda}, \boldsymbol{\mu}) \star \Delta_3^-(\boldsymbol{\lambda}', \boldsymbol{\mu}') = \Delta_6\{(-5d, -3d, -d, d, 3d, 5d), \boldsymbol{\tau}\}$$

is zero-symmetric, where $\boldsymbol{\tau} \propto (r^2, r^3 + r^2, r^3 + r^2 - 1, r^3 + r^2 - 1, r^3 + r^2, r^3 + r^2, r^2)$.

CASE 2. For any real numbers $c$, $d$ and $r$ such that $d > 0$ and $r > 1$, let

(A.4) $\quad \boldsymbol{\mu} = (c, c + 3d, c + 4d) \quad \text{and} \quad \boldsymbol{\lambda} \propto (r\sqrt{r}, (r-1)\sqrt{r+1}, \sqrt{r}),$

(A.5) $\quad \boldsymbol{\mu}' = (c + d, c + 2d, c + 3d) \quad \text{and} \quad \boldsymbol{\lambda}' \propto (r, \sqrt{r + r^2}, 1).$

Then

$$\Delta_3(\boldsymbol{\lambda}, \boldsymbol{\mu}) \star \Delta_3^-(\boldsymbol{\lambda}', \boldsymbol{\mu}') = \Delta_7\{(-3d, -2d, -d, 0, d, 2d, 3d), \boldsymbol{\tau}\}$$

is zero-symmetric, where

$$\boldsymbol{\tau} \propto (r\sqrt{r}, r^2\sqrt{r+1}, r^2\sqrt{r}, (r-1)\sqrt{r+1}, r^2\sqrt{r}, r^2\sqrt{r+1}, r\sqrt{r}).$$

THEOREM A.1. *Let $A$ be the set of all $(\boldsymbol{\lambda}, \boldsymbol{\mu}) \in \Omega_3$ that satisfy any one of the four conditions* (A.2), (A.3), (A.4) *or* (A.5) *for some real numbers $c$, $d$ and $r$ with $d > 0$ and $r > 1$. Let $A^-$ be the set of all $(\boldsymbol{\lambda}, \boldsymbol{\mu}) \in \Omega_3$ such that $((\lambda_3, \lambda_2, \lambda_1), (-\mu_3, -\mu_2, -\mu_1)) \in A$. Finally, let $B$ be the set of all $(\boldsymbol{\lambda}, \boldsymbol{\mu}) \in \Omega_3$ such that $\Delta_3(\boldsymbol{\lambda}, \boldsymbol{\mu})$ is symmetric or $\lambda_1 \lambda_2 \lambda_3 = 0$. Then $\Omega_3^* = \Omega_3 \setminus (A \cup A^- \cup B)$.*

PROOF OF THEOREM A.1. The fact that $\Omega_3^* \subset \Omega_3 \setminus (A \cup A^- \cup B)$ is immediate. Now, let $Y \sim \Delta_3(\boldsymbol{\lambda}, \boldsymbol{\mu})$ with $(\boldsymbol{\lambda}, \boldsymbol{\mu}) \in \Omega_3 \setminus (A \cup A^- \cup B)$ and suppose that $Y - Y'$ is zero-symmetric, where $Y'$ is independent of $Y$ and $Y' \sim \Delta_3(\boldsymbol{\lambda}', \boldsymbol{\mu}')$ for $(\boldsymbol{\lambda}', \boldsymbol{\mu}') \in \Omega_3$. We wish to show that $Y$ and $Y'$ have the same distribution.

Since $Y$ is not symmetric, $Y'$ cannot be a point mass, so we may assume without loss of generality that $\lambda_1'$ and $\lambda_2'$ are positive. To speed things along a bit, we introduce several new variables. Let $\eta_1 < \eta_2 < \cdots < \eta_m$ denote the support points of $Y - Y'$. Furthermore, let $\delta_j = \mu_{j+1} - \mu_j$ and $\delta_j' = \mu_{j+1}' - \mu_j'$, $1 \le j \le 2$, denote the gaps between the elements of $\boldsymbol{\mu}$ and $\boldsymbol{\mu}'$. Finally, let $\alpha_j$ and $\alpha_j'$ equal $\lambda_j/\lambda_1$ and $\lambda_j'/\lambda_1'$, respectively, for $1 \le j \le 3$, so that we may simplify calculations by working with the unnormalized vectors $(1, \alpha_2, \alpha_3)$ and $(1, \alpha_2', \alpha_3')$.

By the zero-symmetry of $Y - Y'$,

$$(A.6) \qquad P(Y - Y' = \eta_j) = P(Y - Y' = \eta_{m-j+1})$$

for $1 \le j \le m$. We now consider the cases $\alpha_3' \ne 0$ and $\alpha_3' = 0$ separately.

CASE A. $\alpha_3' = 0$. By equation (A.6) with $j = 1$, we obtain $\alpha_2' = \alpha_3$. Because $\eta_2 - \eta_1 = \eta_m - \eta_{m-1}$, we have

$$(A.7) \qquad \min\{\delta_1, \delta_1'\} = \min\{\delta_2, \delta_1'\}.$$

If $\delta_1$ and $\delta_2$ are both less than $\delta_1'$, then they must be equal by equation (A.7); if they are both greater than $\delta_1'$, then they must be equal because $\eta_1 + \delta_1$ is the opposite of $-\eta_1 - \delta_2$ by the zero-symmetry of $Y - Y'$. In either case, equation (A.6) with $j = 2$ gives $\alpha_3 = 1$, which is a contradiction because $Y$ cannot be symmetric.

If $\delta_1 = \delta_2 = \delta_1'$, then equation (A.6) with $j = 2$ gives $\alpha_2 + \alpha_3^2 = 1 + \alpha_2 \alpha_3$. Since $\alpha_3 \ne 1$ because $Y$ is not symmetric, this implies that $\alpha_2 = 1 + \alpha_3$ and so either $Y$ or $-Y$ must satisfy condition (A.3), a contradiction.

If $\delta_1 > \delta_1' = \delta_2$, then $\eta_1 + \delta_1$ must be zero because there is no other way that $-(\eta_1 + \delta_1)$ could be attained by $Y - Y'$. This implies that $\delta_1 = 2\delta_2$. From equation (A.6) with $j = 2$, we obtain $\alpha_3^2 + \alpha_2 = 1$ so that $\boldsymbol{\lambda} \propto (1, 1 - \alpha_3^2, \alpha_3)$. Setting $r = 1/\alpha_3$, we see this implies that $Y$ satisfies condition (A.2), a contradiction. The case $\delta_2 > \delta_1' = \delta_1$ leads to a similar contradiction in which $-Y$ satisfies condition (A.2).

By equation (A.7), we have now exhausted case A.

CASE B.   $\alpha_3' \neq 0$. By equation (A.6) with $j = 1$, we have $\alpha_3' = \alpha_3$. In analogy with equation (A.7), we obtain

$$(A.8) \qquad \min\{\delta_1, \delta_2'\} = \min\{\delta_2, \delta_1'\}.$$

We may group the possibilities into the following four categories according to the relative values of $\delta_1$, $\delta_2$, $\delta_1'$ and $\delta_2'$.

B1: $\delta_1 < \delta_2'$ and $\delta_2 < \delta_1'$, or $\delta_1 > \delta_2'$ and $\delta_2 > \delta_1'$.
B2: $\delta_1 < \delta_2'$ and $\delta_2 > \delta_1'$, or $\delta_1 > \delta_2'$ and $\delta_2 < \delta_1'$.
B3: $\delta_1 = \delta_2' = \delta_2 = \delta_1'$.
B4: $\delta_1 > \delta_2'$ and $\delta_2 = \delta_1'$, $\delta_1 < \delta_2'$ and $\delta_2 = \delta_1'$, $\delta_1 = \delta_2'$ and $\delta_2 < \delta_1'$, or $\delta_1 = \delta_2'$ and $\delta_2 > \delta_1'$.

The details of the following arguments are similar to those of Case A, so we omit many of then. In the case B1, $Y$ may be shown to be symmetric, which is a contradiction. In the case B2, we obtain $\boldsymbol{\alpha}' = \boldsymbol{\alpha}$, $\delta_1 = \delta_1'$ and $\delta_2 = \delta_2'$, which means that $Y$ and $Y'$ have the same distribution. In the case B3, either $\alpha_3 = 1$, which leads to a contradiction since $Y$ is then symmetric, or $\boldsymbol{\alpha} = \boldsymbol{\alpha}'$, in which case $Y$ and $Y'$ have the same distribution.

Case B4 implies that three of $\delta_1$, $\delta_2$, $\delta_1'$ and $\delta_2'$ are equal, while the fourth is at least double this common value. For the sake of illustration, suppose that $\delta_1$ is the large value, so that

$$(A.9) \qquad \delta_1 \geq 2\delta_2 = 2\delta_1' = 2\delta_2'.$$

If equality holds in (A.9), then we may show that $Y$ satisfies condition (A.2) and $Y'$ satisfies condition (A.3). On the other hand, if inequality in (A.9) is strict, then $\delta_1$ must equal $3\delta_2$, in which case $Y$ satisfies condition (A.4) and $Y'$ satisfies condition (A.5). Either outcome gives $(\boldsymbol{\lambda}, \boldsymbol{\mu}) \in A$, a contradiction. A similar contradiction occurs when the role of $\delta_1$ is interchanged with that of $\delta_2$, $\delta_1'$ or $\delta_2'$ in (A.9).

Since B1–B4 exhaust Case B, and Case A always leads to a contradiction, we conclude that $Y$ and $Y'$ must have the same distribution.   □

## APPENDIX B: CONSISTENCY PROOFS

PROOF OF LEMMA 1.   For the case of finite $p$, we will show that $\{d(\hat{\boldsymbol{\lambda}}, \hat{\boldsymbol{\mu}})\}^p \to 0$ almost surely. It suffices to prove that

$$\sup_{(\boldsymbol{\lambda}, \boldsymbol{\mu}) \in \Omega_k} |\{d(\boldsymbol{\lambda}, \boldsymbol{\mu})\}^p - \{d_n(\boldsymbol{\lambda}, \boldsymbol{\mu})\}^p| \to 0 \qquad \text{almost surely,}$$

since $d(\boldsymbol{\lambda}^0, \boldsymbol{\mu}^0) = 0$ and $d_n(\hat{\boldsymbol{\lambda}}, \hat{\boldsymbol{\mu}}) \leq d_n(\boldsymbol{\lambda}^0, \boldsymbol{\mu}^0)$ imply that

$$\{d(\hat{\boldsymbol{\lambda}}, \hat{\boldsymbol{\mu}})\}^p \leq \{d(\hat{\boldsymbol{\lambda}}, \hat{\mu})\}^p - \{d_n(\hat{\boldsymbol{\lambda}}, \hat{\mu})\}^p + \{d_n(\boldsymbol{\lambda}^0, \boldsymbol{\mu}^0)\}^p$$

$$\leq |\{d(\hat{\boldsymbol{\lambda}}, \hat{\boldsymbol{\mu}})\}^p - \{d_n(\hat{\boldsymbol{\lambda}}, \hat{\mu})\}^p| + |\{d(\boldsymbol{\lambda}^0, \boldsymbol{\mu}^0)\}^p - \{d_n(\boldsymbol{\lambda}^0, \boldsymbol{\mu}^0)\}^p|$$

$$\leq 2 \sup_{(\boldsymbol{\lambda}, \boldsymbol{\mu}) \in \Omega_k} |\{d(\boldsymbol{\lambda}, \boldsymbol{\mu})\}^p - \{d_n(\boldsymbol{\lambda}, \boldsymbol{\mu})\}^p|.$$

Define the functions

$$\alpha(t) = \sum_{j=1}^{k} \lambda_j \{1 - F_0(\mu_j - t) - F_0(\mu_j + t)\}$$

and

$$\alpha_n(t) = \sum_{j=1}^{k} \lambda_j \{1 - \hat{F}_n(\mu_j - t) - \hat{F}_n(\mu_j + t)\}.$$

Since $x \mapsto x^p - px$ is nonincreasing on $[0, 1]$ for $p \geq 1$, $b^p - a^p \leq p(b - a)$ whenever $0 \leq a \leq b \leq 1$. Therefore,

$$(B.1) \qquad ||a|^p - |b|^p| \leq p|a - b| \qquad \text{for } |a|, |b| \leq 1.$$

Since $|\alpha(t)|$ and $|\alpha_n(t)|$ are both less than or equal to 1, we obtain

$$
\begin{aligned}
|\{d(\boldsymbol{\lambda}, \boldsymbol{\mu})\}^p - \{d_n(\boldsymbol{\lambda}, \boldsymbol{\mu})\}^p| &= \left| \int_{-\infty}^{\infty} (|\alpha(t)|^p - |\alpha_n(t)|^p)\, dt \right| \\
&\leq p \int_{-\infty}^{\infty} |\alpha(t) - \alpha_n(t)|\, dt \\
(B.2) \qquad &\leq p \sum_{j=1}^{k} \lambda_j \int_{-\infty}^{\infty} |F_0(\mu_j - t) - \hat{F}_n(\mu_j - t)|\, dt \\
&\quad + p \sum_{j=1}^{k} \lambda_j \int_{-\infty}^{\infty} |F_0(\mu_j + t) - \hat{F}_n(\mu_j + t)|\, dt \\
&= 2p \int_{-\infty}^{\infty} |F_0(t) - \hat{F}_n(t)|\, dt.
\end{aligned}
$$

Note that taking the supremum over $(\boldsymbol{\lambda}, \boldsymbol{\mu}) \in \Omega_k$ is now irrelevant. Let $g_n(t) = \text{sign}(t)\{\hat{F}_n(t) - F_0(t)\}$ and let $g_n^+(t) = g_n(t)I\{g_n(t) \geq 0\}$ denote its positive part. Let

$$g(t) = \begin{cases} F_0(t), & \text{if } t < 0, \\ 1 - F_0(t), & \text{if } t \geq 0 \end{cases}$$

and note that $0 \leq g_n^+(t) \leq g(t)$. Since $g(t)$ is integrable by the assumption of a finite first moment and $g_n(t) \to 0$ almost surely by the Glivenko–Cantelli theorem, $\int g_n^+(t)\, dt \to 0$ almost surely by the dominated convergence theorem. Furthermore, $\int g_n(t)\, dt = \text{E}_{F_0}|X| - \frac{1}{n}\sum_i |X_i| \to 0$ almost surely by the strong law of large numbers. Since $|g_n(t)| = 2g_n^+(t) - g_n(t)$, this proves that the right-hand side of inequality (B.2) goes to 0 almost surely.

The case of $p = \infty$ is much simpler. Repeated use of the triangle inequality shows that

$$\sup_t |\alpha_n(t)| \leq 2\sup_t |\hat{F}_n(t) - F_0(t)| + \sup_t |\alpha(t)|$$

and the same inequality holds if $\alpha_n(t)$ and $\alpha(t)$ switch places. Therefore,

$$|d(\boldsymbol{\lambda}, \boldsymbol{\mu}) - d_n(\boldsymbol{\lambda}, \boldsymbol{\mu})| \leq 2 \sup_t |\hat{F}_n(t) - F_0(t)|$$

and the Glivenko–Cantelli theorem implies that as $n \to \infty$,

$$\sup_{(\boldsymbol{\lambda}, \boldsymbol{\mu}) \in \Omega_k} |\{d(\boldsymbol{\lambda}, \boldsymbol{\mu})\} - \{d_n(\boldsymbol{\lambda}, \boldsymbol{\mu})\}| \to 0 \qquad \text{almost surely.}$$

Note that no finite first moment assumption is required for $p = \infty$. $\quad \square$

In order to prove Lemma 2, we first define a new function $h(\boldsymbol{\lambda}, \boldsymbol{\mu})$ and show that it is uniformly continuous. The introduction of this function may seem mysterious, but it is designed specifically to resemble the function $d(\boldsymbol{\lambda}, \boldsymbol{\mu})$, while at the same time possessing the crucial uniform continuity property, the importance of which will be discussed further in the proof of Lemma 2.

LEMMA B.1.   *For $1 \leq p < \infty$, the function*

$$
\begin{aligned}
h(\boldsymbol{\lambda}, \boldsymbol{\mu}) = \int_{-\infty}^{\infty} &\left| \sum_{j=1}^{k} \lambda_j \{1 - F_0(\mu_j - t) - F_0(\mu_j + t)\} \right|^p \\
&\times \sum_{j=1}^{k} (e^{-|\mu_j - t|} + e^{-|\mu_j + t|}) \, dt
\end{aligned}
$$
(B.3)

*is uniformly continuous if $G_0$ has finite first moment.*

PROOF.   By inequality (B.1),

$$
\begin{aligned}
|h(\boldsymbol{\lambda}, \boldsymbol{\mu}) - h(\boldsymbol{\lambda}', \boldsymbol{\mu})| &\leq p \sum_{j=1}^{k} |\lambda_j - \lambda_j'| \int_{-\infty}^{\infty} \sum_{j=1}^{k} (e^{-|\mu_j - t|} + e^{-|\mu_j + t|}) \, dt \\
&= 4kp \sum_{j=1}^{k} |\lambda_j - \lambda_j'|.
\end{aligned}
$$

Thus, $h(\boldsymbol{\lambda}, \boldsymbol{\mu})$ is uniformly continuous in $\boldsymbol{\lambda}$. Furthermore, $|a^p c - b^p d| \leq |d(a^p - b^p)| + |a^p(c - d)|$ and inequality (B.1) together imply that $|h(\boldsymbol{\lambda}, \boldsymbol{\mu}) - h(\boldsymbol{\lambda}, \boldsymbol{\mu}')|$ is bounded above by

$$2kp \sum_{j=1}^{k} \lambda_j \int_{-\infty}^{\infty} \{|F_0(\mu_j - t) - F_0(\mu_j' - t)| + |F_0(\mu_j + t) - F_0(\mu_j' + t)|\} \, dt$$

$$+ \sum_{j=1}^{k} \int_{-\infty}^{\infty} (|e^{-|\mu_j - t|} - e^{-|\mu_j' - t|}| + |e^{-|\mu_j + t|} - e^{-|\mu_j' + t|}|) \, dt.$$

Since each of the above integrals is over the whole real line, each depends on $\mu_j$ and $\mu'_j$ solely through the difference $\mu_j - \mu'_j$. The dominated convergence theorem implies that each integral tends to zero as $\boldsymbol{\mu} - \boldsymbol{\mu}' \to \mathbf{0}$, which establishes that $h(\boldsymbol{\lambda}, \boldsymbol{\mu})$ is also uniformly continuous in $\boldsymbol{\mu}$. $\quad\square$

PROOF OF LEMMA 2. We first consider the case of finite $p$ (the case $p = \infty$ is much easier). Since $h(\boldsymbol{\lambda}, \boldsymbol{\mu})$ defined in equation (B.3) is bounded above by $2k\{d(\boldsymbol{\lambda}, \boldsymbol{\mu})\}^p$, if Lemma 2 is false, then there exist $\varepsilon > 0$ and a sequence

$$\{(\boldsymbol{\lambda}^n, \boldsymbol{\mu}^n)\}_{n=1}^\infty \subset \{(\boldsymbol{\lambda}, \boldsymbol{\mu}) \in \Omega_k : \|(\boldsymbol{\lambda}, \boldsymbol{\mu}) - (\boldsymbol{\lambda}^0, \boldsymbol{\mu}^0)\| > \varepsilon\}$$

such that $h(\boldsymbol{\lambda}^n, \boldsymbol{\mu}^n) \to 0$ as $n \to \infty$. We now show that this leads to a contradiction.

Passing to a subsequence if necessary, we assume without loss of generality that $\boldsymbol{\lambda}^n \to \boldsymbol{\lambda}^*$ and that each of the sequences $\mu_j^n$ has a limit, either finite or infinite. By the uniform continuity of $h(\boldsymbol{\lambda}, \boldsymbol{\mu})$ (Lemma B.1), we have

$$(B.4) \qquad\qquad h(\boldsymbol{\lambda}^*, \boldsymbol{\mu}^n) \to 0.$$

Note that we cannot obtain an analogous expression using $d(\boldsymbol{\lambda}, \boldsymbol{\mu})$ instead of $h(\boldsymbol{\lambda}, \boldsymbol{\mu})$ because $d$ is not uniformly continuous.

A standard result in analysis states that if $f_n \to f$ in $L^1(R)$, then $f_n$ has a subsequence that converges to $f$ almost everywhere (a.e.) with respect to Lebesgue measure. Thus, passing to a subsequence if necessary, we see that (B.4) implies

$$(B.5) \qquad \left| \sum_{j=1}^k \lambda_j^* \{1 - F_0(\mu_j^n - t) - F_0(\mu_j^n + t)\} \right|^p \sum_{j=1}^k (e^{-|\mu_j^n - t|} + e^{-|\mu_j^n + t|}) \to 0 \text{ a.e.}$$

If some of the sequences $\{\mu_j^n\}$ are bounded, say $\mu_j^n \to \mu_j^*$, replacing these sequences by their limits does not change (B.4) because of the uniform continuity of $h(\boldsymbol{\lambda}, \boldsymbol{\mu})$. Furthermore, in this case, the second sum in (B.5) is bounded away from zero, which implies that the expression inside the absolute value symbols tends to zero for almost all $t$. We conclude that if $\{j : \mu_j^n \to \mu_j^*\}$ is nonempty, then

$$(B.6) \qquad \sum_{j : \mu_j^n \to \mu_j^*} \lambda_j^* \{1 - F_0(\mu_j^* - t) - F_0(\mu_j^* + t)\} + \sum_{j : \mu_j^n \to -\infty} \lambda_j^* - \sum_{j : \mu_j^n \to \infty} \lambda_j^* = 0 \text{ a.e.}$$

Letting $t \to \infty$ in equation (B.6) gives

$$\sum_{j : \mu_j^n \to -\infty} \lambda_j^* - \sum_{j : \mu_j^n \to \infty} \lambda_j^* = 0.$$

Hence, equation (B.6) implies that

$$\sum_{j\,:\,\mu_j^n \to \mu_j^*} \lambda_j^* \{1 - F_0(\mu_j^* - t) - F_0(\mu_j^* + t)\} = 0 \text{ a.e.},$$

which contradicts the assumption that $F_0$ is $k$-identifiable (note that we may assume without loss of generality that $\lambda_j^* \neq 0$ whenever $\mu_j^n \to \mu_j^*$; otherwise, this $j$th component may be entirely ignored). Therefore, every sequence $\{\mu_j^n\}$ goes to $\pm\infty$ as $n \to \infty$.

Next, take $C$ to be an arbitrary constant that is not contained in the set $\{\mu_1^0, \ldots, \mu_k^0\}$. Fix $j_0$ such that $\lambda_{j_0}^* \neq 0$, then define $a_j^n = \mu_j^n - \mu_{j_0}^n + C$ and $b_j^n = \mu_j^n + \mu_{j_0}^n - C$ for $1 \leq j \leq k$ and $n \geq 1$. Under the change of variable $s = C - \mu_{j_0}^n - t$,

$$h(\boldsymbol{\lambda}^*, \boldsymbol{\mu}^n) = \int_{-\infty}^{\infty} \left| \sum_{j=1}^{k} \lambda_j^* \{1 - F_0(a_j^n - s) - F_0(b_j^n + s)\} \right|^p$$

$$\times \sum_{j=1}^{k} (e^{-|a_j^n - s|} + e^{-|b_j^n + s|}) \, ds.$$

Thus, passing to a subsequence if necessary, the argument leading to (B.5) implies that

$$(\text{B.7}) \quad \left| \sum_{j=1}^{k} \lambda_j^* \{1 - F_0(a_j^n - s) - F_0(b_j^n + s)\} \right|^p \sum_{j=1}^{k} (e^{-|a_j^n - s|} + e^{-|b_j^n + s|}) \to 0 \text{ a.e.}$$

Since $a_{j_0}^n = C$ for all $n$ by definition, the second sum in (B.7) is bounded away from zero, implying that

$$(\text{B.8}) \quad \sum_{j=1}^{k} \lambda_j^* \{1 - F_0(a_j^n - s)\} - \sum_{j=1}^{k} \lambda_j^* F_0(b_j^n + s) \to 0 \text{ a.e.}$$

Passing to a subsequence if necessary, we assume that all of the sequences $a_j^n$ and $b_j^n$ have limits, either finite or infinite. For any of these sequences whose limit is finite, we denote this limit by $a_j^*$ or $b_j^*$. Thus, we may decompose the sum in (B.8) into parts according to the limits of the $a_j^n$ and $b_j^n$, as follows:

$$(\text{B.9}) \quad \sum_{j\,:\,a_j^n \to a_j^*} \lambda_j^* \{1 - F_0(a_j^* - s)\} - \sum_{j\,:\,b_j^n \to b_j^*} \lambda_j^* F_0(b_j^* + s) + \sum_{j\,:\,a_j^n \to -\infty} \lambda_j^*$$

$$- \sum_{j\,:\,b_j^n \to \infty} \lambda_j^* = 0 \text{ a.e.}$$

Letting $s \to -\infty$ in equation (B.9) gives

$$\sum_{j\,:\,a_j^n \to -\infty} \lambda_j^* - \sum_{j\,:\,b_j^n \to \infty} \lambda_j^* = 0.$$

Therefore, equation (B.9) implies that

$$(B.10) \qquad \sum_{j: a_j^n \to a_j^*} \lambda_j^* \{1 - F_0(a_j^* - s)\} = \sum_{j: b_j^n \to b_j^*} \lambda_j^* F_0(b_j^* + s) \text{ a.e.}$$

Let $k_1$ and $k_2$ denote the number of finite $a_j^*$ and $b_j^*$, respectively. Note that $k_1$ and $k_2$ must be positive by equation (B.10) because $a_{j_0}^* = C$ and $\lambda_{j_0}^* \neq 0$. Furthermore, since $|\mu_j^n| \to \infty$ for each $j$, at most one of $\{a_j^n\}$ and $\{b_j^n\}$ can remain bounded; thus, $k_1 + k_2 \leq k$. (This gives an immediate contradiction if $k = 1$.)

Equation (B.10) is an identity of distribution functions. It states that

$$(B.11) \qquad F_0^-(t) \star \Delta_{k_1}(\boldsymbol{\alpha}, \boldsymbol{\xi}) = F_0(t) \star \Delta_{k_2}(\boldsymbol{\beta}, \boldsymbol{\eta}),$$

where $\Delta_{k_1}(\boldsymbol{\alpha}, \boldsymbol{\xi})$ is the distribution function supported on the finite $a_j^*$ with weights proportional to the corresponding $\lambda_j^*$ and $\Delta_{k_2}(\boldsymbol{\beta}, \boldsymbol{\eta})$ is the distribution function supported on the finite $-b_j^*$ with weights proportional to the corresponding $\lambda_j^*$. Recall that $F_0 = G_0 \star \Delta_k(\boldsymbol{\lambda}^0, \boldsymbol{\mu}^0)$. Define independent random variables $Z \sim G_0$, $Y \sim \Delta_k(\boldsymbol{\lambda}^0, \boldsymbol{\mu}^0)$, $W_1 \sim \Delta_{k_1}(\boldsymbol{\alpha}, \boldsymbol{\xi})$ and $W_2 \sim \Delta_{k_2}(\boldsymbol{\beta}, \boldsymbol{\eta})$. Then in terms of characteristic functions, equation (B.11) becomes $\phi_{W_1} \phi_{-Y} \phi_Z = \phi_{W_2} \phi_Y \phi_Z$. We may cancel $\phi_Z$ from both sides because $\phi_{W_1} \phi_{-Y}$ and $\phi_{W_2} \phi_Y$ are analytic functions that agree whenever $\phi_Z \neq 0$. This leads to

$$(B.12) \qquad \frac{1}{2}(\phi_{-W_1} \phi_Y + \phi_{W_1} \phi_{-Y}) = \frac{1}{2}(\phi_{-W_1} + \phi_{W_2}) \phi_Y.$$

The left-hand side of equation (B.12) is a real function, which means that

$$\frac{1}{2}\{\Delta_{k_1}^-(\boldsymbol{\alpha}, \boldsymbol{\xi}) + \Delta_{k_2}(\boldsymbol{\beta}, \boldsymbol{\eta})\} \star \Delta_k(\boldsymbol{\lambda}^0, \boldsymbol{\mu}^0)$$

is zero-symmetric. Since $k_1 + k_2 \leq k$ and $(\boldsymbol{\lambda}^0, \boldsymbol{\mu}^0) \in \Omega_k^*$, we have

$$(B.13) \qquad \frac{1}{2}\{\Delta_{k_1}^-(\boldsymbol{\alpha}, \boldsymbol{\xi}) + \Delta_{k_2}(\boldsymbol{\beta}, \boldsymbol{\eta})\} = \Delta_k^-(\boldsymbol{\lambda}^0, \boldsymbol{\mu}^0)$$

by Theorem 1. But this is impossible, because the distribution on the left side of equation (B.13) assigns nonzero weight to the point $-a_1 = -C$ and $C$ was chosen specifically so that the distribution on the right assigns no weight at $-C$. This completes the proof for finite $p$.

If $p = \infty$, then

$$d(\boldsymbol{\lambda}, \boldsymbol{\mu}) = \sup_t \left| \sum_{j=1}^k \lambda_j \{1 - F_0(\mu_j - t) - F_0(\mu_j + t)\} \right|$$

is uniformly continuous in $\boldsymbol{\lambda}$, so the $h(\boldsymbol{\lambda}, \boldsymbol{\mu})$ function is unnecessary. Thus, we assume the existence of a sequence $(\boldsymbol{\lambda}^n, \boldsymbol{\mu}^n)$ such that $\boldsymbol{\lambda}^n \to \boldsymbol{\lambda}^*$, each $\boldsymbol{\mu}_j^n$ has a limit and $d(\boldsymbol{\lambda}^*, \boldsymbol{\mu}^n) \to 0$. Equation (B.6) is true because $F_0(t)$ is almost everywhere continuous; equation (B.8) is trivially true. The rest of the proof for $p = \infty$ is identical to the proof for $p < \infty$, and we arrive at a contradiction.  $\square$

D. R. Hunter
T. P. Hettmansperger
Department of Statistics
Pennsylvania State University
University Park, Pennsylvania 16802-2111
USA
E-mail: dhunter@stat.psu.edu
        tph@stat.psu.edu

S. Wang
School of Public Health
Yale University
New Haven, Connecticut 06520
USA
E-mail: swang@masal.med.yale.edu